\documentclass[11pt,reqno]{amsart}

\usepackage[alphabetic,nobysame]{amsrefs}
\usepackage[font=footnotesize]{caption}
\usepackage{bm}

\usepackage[margin=1.32in]{geometry}
\usepackage[T1]{fontenc}
\usepackage[utf8]{inputenc}
\usepackage{lmodern}
\usepackage{microtype}

\usepackage{amsmath,amsfonts,amssymb,amsthm,mathtools}

\usepackage[dvipsnames,svgnames]{xcolor}
\usepackage{tikz}
\usepackage{hyperref}
\usepackage{enumitem}
\usepackage[nameinlink,capitalize,noabbrev]{cleveref}

\hypersetup{
  colorlinks=true,
  linkcolor=MidnightBlue,
  citecolor=ForestGreen,
  urlcolor=RoyalBlue
}

\theoremstyle{plain}
\newtheorem{theorem}{Theorem}[section]
\newtheorem{proposition}[theorem]{Proposition}
\newtheorem{lemma}[theorem]{Lemma}
\theoremstyle{definition}
\newtheorem{definition}[theorem]{Definition}
\theoremstyle{remark}
\newtheorem{remark}[theorem]{Remark}
\newtheorem{example}[theorem]{Example}

\theoremstyle{plain}
\newtheorem{introtheorem}{Theorem}
\crefname{introtheorem}{theorem}{theorems}
\Crefname{introtheorem}{Theorem}{Theorems}

\AddToHook{env/lemma/begin}{\crefalias{theorem}{lemma}}
\AddToHook{env/proposition/begin}{\crefalias{theorem}{proposition}}
\AddToHook{env/definition/begin}{\crefalias{theorem}{definition}}
\AddToHook{env/remark/begin}{\crefalias{theorem}{remark}}
\AddToHook{env/example/begin}{\crefalias{theorem}{example}}

\newcommand{\RR}{\mathbb{R}}
\newcommand{\EE}{\mathbb{E}}
\newcommand{\PP}{\mathbb{P}}
\newcommand{\eps}{\varepsilon}

\DeclareMathOperator{\Hess}{Hess}
\DeclareMathOperator{\Cut}{Cut}
\DeclareMathOperator{\Ric}{Ric}
\DeclareMathOperator{\Len}{Len}

\newcommand{\cV}{\mathcal{V}}
\newcommand{\cA}{\mathcal{A}}
\newcommand{\cF}{\mathcal{F}}

\title[Free-time minimizers for $N$-body problems and a splitting theorem]{Absence of free-time minimizers for\\ $\bm{N}$-body problems and a splitting theorem}

\author[R. Assouline]{Rotem Assouline}
\address{Rotem Assouline\newline\indent Massachusetts Institute of Technology, Department of Mathematics\newline\indent Cambridge, MA 02139-4307, USA}
\email{rotem\_a@mit.edu}

\author[M. Mazzucchelli]{Marco Mazzucchelli}
\address{Marco Mazzucchelli\newline\indent Sorbonne Université, Université Paris Cité, CNRS, IMJ-PRG\newline\indent F-75005 Paris, France\newline\indent\indent and\newline\indent Institute for Advanced Study, School of Mathematics\newline\indent  1 Einstein Drive, Princeton, NJ 08540, USA}
\email{mmazzucchelli@ias.edu}

\date{September 24, 2026}

\keywords{Lagrangian systems, free-time minimizers, $N$-body problem, splitting theorem}

\subjclass[2020]{70F10, 53C20}

\begin{document}

\begin{abstract}
  For the $N$-body problem in dimension $d \ge 3$, we give an alternative proof of a theorem of da Luz and Maderna on the nonexistence of free-time minimizers defined on the entire real line, and extend it from the Newtonian potential to a larger class of singular subharmonic pairwise interaction potentials.
  The proof is based on the following generalization of the Cheeger--Gromoll splitting theorem: let $(M,g)$ be a connected Riemannian manifold and let $U$ be a smooth positive function such that $\Ric_g \ge \tfrac{\Delta U}{2U}g$, and the metrics $Ug$ and $U^{-1}g$ are both complete; if the metric $Ug$ admits a line, then the function $U$ is constant and $(M,g)$ splits off that line. We apply this splitting theorem to the Jacobi--Maupertuis metric associated with the $N$-body problem;
  a quantitative version of Marchal's theorem enables us to reduce to the case of smooth potentials.
\end{abstract}

\maketitle

\section{Introduction}

Consider the $N$-body problem in Euclidean space $\RR^d$, with $N \ge 2$, $d \ge 3$, masses $m_1,\dots,m_N > 0$, and potential $V : (\RR^d)^N \to [-\infty,0)$;
a natural choice is the Newtonian potential
$V(x):=-\sum m_im_j\|x_i - x_j\|_{\RR^d}^{-1}$.

A curve $\gamma = (\gamma_1,\dots,\gamma_N):I \to (\RR^d)^N$ defined on an interval $I\subseteq \RR$ is said to be a \emph{free-time minimizer} if, for every $t,t'\in I$ with $t<t'$, the action
\[
  \int_t^{t'}\left(\frac12\sum_{i=1}^Nm_i\|\dot\gamma_i(s)\|_{\RR^d}^2 - V(\gamma(s))\right)ds
\]
is minimal among absolutely continuous curves defined on any finite interval and joining $\gamma(t)$ to $\gamma(t')$. 
A free-time minimizer is \emph{entire} if it is defined on all of $\RR$.

In Section~\ref{sec:V}, we define a class $\cV_{N,d}$ of potentials which contains the Newtonian one. More generally, $\cV_{N,d}$ contains  every potential of the form 
\[-\sum_{i<j}c_{ij}\,\|x_i - x_j\|_{\RR^d}^{-\alpha_{ij}}\] 
with $c_{ij} > 0$, $\alpha_{ij}\in(0,1]$ if $d=3$, and $\alpha_{ij}\in(0,2)$ if $d>3$. It also contains potentials in which the individual interactions are not homogeneous; the key property of the potentials in $\cV_{N,d}$ is \emph{subharmonicity}. For this class we prove the following theorem.

\begin{introtheorem}\label{thm:mb}
  There are no entire free-time minimizers for the $N$-body problem in $\RR^d \,(d \ge 3)$ with any potential $V \in \cV_{N,d}$.
\end{introtheorem}

In the Newtonian case, \Cref{thm:mb} is due to da Luz and Maderna \cite[Theorem 2]{DaLuzMaderna2014} for $d \ge 2$. Their argument, which relies on homogeneity of the potential, uses rescaling to show that the existence of an entire free-time minimizer contradicts Marchal's theorem.
The proof we give here is of a different nature: we deduce \Cref{thm:mb} from the following extension of the Cheeger--Gromoll splitting theorem \cite{CheegerGromoll1971,EschenburgHeintze1984,Petersen2016}.
Recall that a \emph{line} on a Riemannian manifold is a unit-speed minimizing geodesic defined on all of $\RR$.

\begin{introtheorem}\label{thm:main}
  Let $(M,g)$ be a connected Riemannian manifold of dimension $n \ge 2$ and let $U : M \to (0,\infty)$ be a smooth function. Suppose that
  \begin{equation}\label{eq:DeltaVcondition}
    \Ric_g \ge \frac{\Delta U}{2U} g,
  \end{equation}
  and that the conformal metrics $U g$ and $U^{-1}g$ are both complete.
  If there exists a line for the metric $Ug$, then the function $U$ is constant and $(M,g)$ admits an isometric splitting $M \cong M'\times \RR$, where $M'$ is a complete Riemannian manifold of nonnegative Ricci curvature.
\end{introtheorem}

\begin{remark}
  \Cref{prop:superharmonic} is the only place where completeness of $U^{-1}g$ is used. By examining its proof, one sees that it is sufficient to assume instead that every minimizing geodesic ray of the metric $Ug$ has infinite length with respect to the metric $U^{-1}g$.
\end{remark}

The proof of \Cref{thm:main} is essentially the same as that of the Cheeger--Gromoll splitting theorem: it is likewise based on analysis of Busemann functions and a Laplacian comparison principle for the distance function of the metric $Ug$. 
However, \Cref{thm:main} cannot be proved by directly applying the Cheeger--Gromoll splitting theorem to the metric $Ug$, since condition \eqref{eq:DeltaVcondition} does not imply that it has nonnegative Ricci curvature. 

We deduce \Cref{thm:mb} from \Cref{thm:main} by taking $g$ to be the mass metric and setting $U=-V$, so that $2Ug$ is the Jacobi-Maupertuis metric.
Free-time minimizers for the $N$-body problem with potential $V$ are then precisely minimizing geodesics of the metric $2Ug$, parametrized by unit speed with respect to the dual metric $(2U)^{-1}g$.

Busemann-type functions and minimizing rays for the $N$-body problem were considered in \cite{Contreras2001,Maderna2013,PercinoSanchezMorgado2014,MadernaVenturelli2020,MadernaVenturelli2026,BurgosMaderna2021,Montgomery2024}. Lee \cite{Lee2018} employed weak KAM theory to obtain almost-splitting theorems for mechanical Lagrangians.  Comparison theorems for the Laplacians of solutions to the Hamilton--Jacobi equation for Tonelli Lagrangians were established in~\cite{Ohta2014,Assouline2024}.

In recent works, condition \eqref{eq:DeltaVcondition} has been called a \emph{spectral Ricci curvature bound}. For splitting theorems under this condition when $g$ itself is complete, see \cite{AntonelliPozzettaXu2024,HongWang2026} and references therein.

\subsection*{Organization of the paper}
In \Cref{sec:V} we deduce \Cref{thm:mb} from \Cref{thm:main} by constructing regularized potentials for the $N$-body problem with potential $V \in \cV_{N,d}$ and verifying the hypotheses of \Cref{thm:main} for them. In \Cref{sec:proof} we prove \Cref{thm:main}.

\subsection*{Acknowledgements}
The first author is supported by the Fondation Sciences Mathématiques de Paris (FSMP) and the Rothschild Fellowship (Yad Hanadiv). The second author is partially supported by the ANR project QuanQual (ANR-26-CE40-0676), and completed part of this work while in residence at the Institute for Advanced Study during the Special Year on Conformally Symplectic Dynamics and Geometry.

\section[Proof of Theorem 1]{Proof of \Cref{thm:mb}}\label{sec:V}

We begin by introducing the class $\cV_{N,d}$ of potentials mentioned in \Cref{thm:mb}.

For $d \ge 3$, denote by $\cF_d$ the class of smooth functions $f : (0,\infty)\to(0,\infty)$ satisfying the following three conditions.
\begin{enumerate}[label=\textup{(f\arabic*)}, ref=f\arabic*, leftmargin=*]
  \setlength{\itemsep}{5pt}
  \item\label{eq:f-superharmonic} the function $r \mapsto r^{d-1}f'(r)$ is nonincreasing,
  \item\label{eq:f-bounds} there exist $a \in (0,2)$ and $C > 0$ such that
    \[
      C^{-1}(1+r)^{-a} \le f(r) \le C\left(1 + r^{-a}\right), \qquad r > 0 ,
    \]
  \item\label{eq:f-singular} the function $f$ satisfies
    \[
      \lim_{r \searrow 0}\frac{f(r)}{\log(1/r)} = +\infty.
    \]
\end{enumerate}

Condition \eqref{eq:f-superharmonic} implies that $f(\|\cdot\|_{\RR^d})$ is superharmonic on $\RR^d$. By \eqref{eq:f-superharmonic} and \eqref{eq:f-singular}, $f$ is strictly decreasing and diverges to $+\infty$ as $r\searrow 0$; we extend any such $f$ to a continuous function $f:[0,\infty)\to(0,\infty]$ by setting $f(0)=\infty$. When $d=3$, condition \eqref{eq:f-bounds} follows from \eqref{eq:f-superharmonic}.

\begin{definition}
  For $N \ge 2$ and $d \ge 3$, denote by $\cV_{N,d}$ the class of potentials $V : (\RR^d)^N \to [-\infty,0)$ of the form
  \[
    V(x)  =  -\sum_{i<j}f_{ij}\!\left(\|x_i - x_j\|_{\RR^d}\right) , \qquad x = (x_1,\dots,x_N) \in (\RR^d)^N,
  \]
  with $f_{ij} \in \cF_d$ for all $1 \le i < j \le N$.
\end{definition}

\begin{example}[Superpositions of homogeneous pairwise potentials are in $\cV_{N,d}$]
  Let $d \ge 3$ and $N \ge 2$, and let $0 < A < 2$ if $d \ge 4$ and $0 < A \le 1$ if $d=3$. Let $\{\nu_{ij}\}_{i<j}$ be nonzero finite measures on $(0,A]$, and put
  \[
    f_{ij}(r) : = \int_{(0,A]} r^{-\alpha}\,d\nu_{ij}(\alpha) .
  \]
  Then each $f_{ij}$ is smooth and satisfies \eqref{eq:f-superharmonic} and \eqref{eq:f-bounds} with $a = A$, and choosing $\alpha_{0} > 0$ with $\nu_{ij}\!\left([\alpha_{0},A]\right) > 0$ gives $f_{ij}(r) \ge \nu_{ij}\!\left([\alpha_{0},A]\right)\,r^{-\alpha_{0}}$ for $r \le 1$, so \eqref{eq:f-singular} holds as well. 
  Thus the potential
  \[
    V(x) = -\sum_{i<j}\int_{(0,A]} \|x_i-x_j\|_{\RR^d}^{-\alpha}\,d\nu_{ij}(\alpha), \qquad x = (x_1,\dots,x_N) \in (\RR^d)^N
  \]
  belongs to the class $\cV_{N,d}$.
  The Newtonian potential is the case  $\nu_{ij} = m_im_j \cdot \delta_1$, and measures $\nu_{ij}$ supported on two points give the quasihomogeneous potentials introduced in \cite{Diacu1996} (in the special case where the exponents are $\in(0,A]$).
\end{example}

For the remainder of this section, fix $N \ge 2$, $d \ge 3$ and $m_1,\dots,m_N > 0$, and let $\langle\cdot,\cdot\rangle_m$ denote the mass inner product on $(\RR^d)^N$:
\[
  \langle x,y\rangle_m = \sum_{i=1}^Nm_i\,\langle x_i,y_i\rangle_{\RR^d}, \qquad x = (x_1,\dots,x_N),\ y = (y_1,\dots,y_N) \in (\RR^d)^N,
\]
where $\langle\cdot,\cdot\rangle_{\RR^d}$ is the standard inner product on $\RR^d$. Write $\|\cdot\|_m$ for the norm induced by $\langle\cdot,\cdot\rangle_m$.
Fix also a potential $V \in \cV_{N,d}$.

Since \eqref{eq:f-bounds} persists (with a possibly larger $C$) when the exponent $a$ is increased, we may choose a single $a \in (0,2)$ and $C > 0$ for which \eqref{eq:f-bounds} is satisfied by $f_{ij}$ for all $i<j$.

In order to prove \Cref{thm:mb}, we would like to apply \Cref{thm:main}, taking $g$ to be the mass metric 
\[
  g_m : = \left<\cdot,\cdot\right>_m
\]
and the conformal factor to be
\[
  U : = - V;
\]
the metric
\[
  h : = 2Ug_m
\]
is then the \emph{Jacobi--Maupertuis metric}.
However, the conformal factor $2U$ is singular at collisions (that is, at those points $x\in(\RR^d)^N$ such that $x_i=x_j$ for some $i<j$),
and therefore the metric $h$ is only defined on the complement of the collision set. Nevertheless, we can still define $d_h$ as an extended real-valued function on the whole $(\RR^d)^N\times(\RR^d)^N$ by
\begin{align*}
  d_h(x,y)=\inf_{\gamma}\Len_h(\gamma),
\end{align*}
where the infimum ranges over all absolutely continuous curves $\gamma:[t,t']\to (\RR^d)^N$ joining $x$ and $y$, and the (possibly infinite) length $\Len_h(\gamma)$ is defined as usual as the integral of the singular norm $\|\cdot\|_h$ along $\gamma$. Here, if $x\in(\RR^d)^N$ is a collision and $v\in T_x(\RR^d)^N$, we set $\|v\|_h:=+\infty$ if $v\neq0$ and $\|v\|_h=0$ if $v=0$.
In \Cref{lem:mane} we prove that, in fact, the distance function $d_h$ is finite everywhere.

We will construct a family $V_\eps$ of smooth, finite potentials approximating $V$, verify the hypotheses of \Cref{thm:main} for the conformal factors $U_\eps : = -V_\eps$, and show that an entire free-time minimizer for the  potential $V$ would also be an entire free-time minimizer for the potential $V_\eps$, if $\eps$ is sufficiently small. We will then apply \Cref{thm:main} to $g_m$ and $U_\eps$ and obtain a contradiction.

For $i<j$ and $\eps \in (0,1]$, let $f^{\eps}_{ij} : [0,\infty) \to (0,\infty)$ be a nonincreasing function such that $f^{\eps}_{ij}(\|\cdot\|_{\RR^d})$ is smooth and superharmonic on $\RR^d$ and such that
\begin{equation}\label{eq:reg}
  \begin{split}
    f^{\eps}_{ij}(r) = f_{ij}(r), &\qquad \forall r \ge \eps,\\ \qquad f_{ij}(\eps) \le f^{\eps}_{ij}(r) \le f_{ij}(r), &\qquad \forall r\in(0,\eps). 
  \end{split}
\end{equation}
Such $f^{\eps}_{ij}$ can be constructed by multiplying $f'_{ij}$ by a smooth nondecreasing cutoff function, and then integrating back.
We define a smooth potential $V_\eps:(\RR^d)^N\to(-\infty,0)$ by 
\[
  V_\eps(x) : = -\sum_{i<j}f^{\eps}_{ij}\!\left(\|x_i - x_j\|_{\RR^d}\right) , \qquad x \in (\RR^d)^N ,
\]
and set
\[
  U_\eps : = -V_\eps, \qquad h_\eps : = 2U_\eps g_m, \qquad h_\eps^* : = (2U_\eps)^{-1} g_m .
\]
For an absolutely continuous curve $\gamma:[t_1,t_2] \to (\RR^d)^N$, we denote by
\[
  \cA_\eps(\gamma) : = \int_{t_1}^{t_2}\left(\tfrac12\left\|\dot\gamma(t)\right\|_m^2 - V_\eps(\gamma(t))\right)dt
\]
its action with respect to the potential $V_\eps$.
Free-time minimizers for the $N$-body problem with potential $V_\eps$ will be called \emph{$\eps$-minimizers}.
We also set 
\[f^{0}_{ij} = f_{ij},\quad V_0 = V,\quad U_0=U=-V,\quad h_0 = h,\quad h_0^*=(2U)^{-1}g_m,\] 
and denote by $\cA_0$ the action corresponding to the potential $V$. The $0$-minimizers (i.e., free-time minimizers for the $N$-body problem with potential $V$) will be called simply \emph{minimizers}.

By the arithmetic--geometric mean inequality, for every $\eps \in [0,1]$ and every absolutely continuous curve $\gamma : [t_1,t_2] \to (\RR^d)^N$, we have
\begin{equation}\label{eq:length-vs-action}
  \cA_\eps(\gamma) \ge \int_{t_1}^{t_2}\!\sqrt{2U_\eps(\gamma(t))}\,\left\|\dot\gamma(t)\right\|_m dt  =  \Len_{h_\eps}(\gamma).
\end{equation}
If $\cA_\eps(\gamma)$ is finite then there is equality in \eqref{eq:length-vs-action} if and only if $\|\dot\gamma\|_m^2 = 2U_\eps(\gamma)$ almost everywhere.
Thus, taking the infimum in \eqref{eq:length-vs-action} over all absolutely continuous curves joining two given points, we obtain the well-known \emph{Maupertuis' principle}:
for every $\eps \in [0,1]$ and all $x,x' \in (\RR^d)^N$, we have
\begin{equation}\label{eq:maupertuis}
  \inf_\eta \cA_\eps(\eta)  =  d_{h_\eps}(x,x'),
\end{equation}
where the infimum ranges over all absolutely continuous curves $\eta:I\to (\RR^d)^N$ defined on any compact interval $I$ and joining $x$ and $x'$.
For $\eps\in(0,1]$, a curve $\gamma : I \to (\RR^d)^N$ is an $\eps$-minimizer if and only if it is a reparametrization with unit $h_\eps^*$-speed of a minimizing $h_\eps$-geodesic. Therefore any such $\gamma$ satisfies
\begin{equation}\label{eq:unit-speed}
  \left\|\dot\gamma(t)\right\|_m^2  =  2U_\eps(\gamma(t)), \qquad \forall  t \in I,
\end{equation}
and
\begin{equation}\label{eq:calibrated}
  \cA_\eps\big(\gamma|_{[t,t']}\big)  =  \Len_{h_\eps}\!\big(\gamma|_{[t,t']}\big)   = d_{h_\eps}(\gamma(t),\gamma(t')), \qquad \forall t,t'\in I\mbox{ with } t < t'.
\end{equation}
For $\eps = 0$, these assertions still hold for minimizers that avoid collisions.

\begin{lemma}\label{reg:hypotheses}
  For every $\eps \in (0,1]$, the hypotheses of \Cref{thm:main} are satisfied with $U = 2U_\eps$ and $g=g_m$.
\end{lemma}

\begin{proof}
  The function $U_\eps$ is smooth and positive. Since the mass metric is flat, condition \eqref{eq:DeltaVcondition} reads $\Delta U_\eps \le 0$, which follows from the fact that each $f^{\eps}_{ij}(\|\cdot\|_{\RR^d})$ is superharmonic since $\Delta_m = \sum_im_i^{-1}\Delta_{x_i}$.
  By \eqref{eq:reg} and \eqref{eq:f-bounds},
  \[
    f^{\eps}_{12}(r) \ge f_{12}(\max\{\eps,r\}) \gtrsim (1+r)^{-a}, \qquad \forall r > 0,
  \]
  and therefore
  \begin{equation}\label{eq:Ueps-lower}
    U_\eps(x) \ge f^{\eps}_{12}(\|x_1 - x_2\|_{\RR^d}) \gtrsim \left(1+\|x\|_m\right)^{-a}, \qquad x \in (\RR^d)^N.
  \end{equation}
  Let $x, x' \in (\RR^d)^N$ satisfy $\|x'\|_m \ge \|x\|_m$ and let $\gamma : [0,T] \to (\RR^d)^N$ be an absolutely continuous curve joining $x$ to $x'$, parametrized by arclength for the mass metric. Then the function $\|\gamma\|_m$ is $1$-Lipschitz, so $\big|\tfrac{d}{ds}\|\gamma\|_m\big| \le 1$ almost everywhere, and its image contains $\left[\|x\|_m,\|x'\|_m\right]$. Therefore, by~\eqref{eq:Ueps-lower},
  \begin{align*}
    \Len_{h_\eps}(\gamma)  & =  \int_0^T\!\sqrt{2U_\eps(\gamma(s))}\,ds \ \gtrsim\ \int_0^T\!\left(1+\|\gamma(s)\|_m\right)^{-a/2}ds \\
    &\ge\ \int_0^T\!\left(1+\|\gamma(s)\|_m\right)^{-a/2}\left|\tfrac{d}{ds}\|\gamma(s)\|_m\right|ds \ge \int_{\|x\|_m}^{\|x'\|_m}\!\left(1+r\right)^{-a/2}dr.
  \end{align*}
  Taking the infimum over $\gamma$ we get 
  \[
    d_{h_\eps}(x,x') \gtrsim \int_{\|x\|_m}^{\|x'\|_m}\!\left(1+r\right)^{-a/2}dr \xrightarrow{\|x'\|_m \to \infty} \infty
  \]
  because $a < 2$. Thus closed $h_\eps$-balls are compact, and therefore the metric $h_\eps$ is complete.
  Finally, each $f^{\eps}_{ij}$ is continuous and nonincreasing on $[0,\infty)$, hence bounded; thus $U_\eps$ is bounded, so completeness of $h_\eps^*$ follows from that of $g_m$.
\end{proof}

While the Jacobi-Maupertuis metric $h$ is singular at collisions, the following lemma implies that the distance $d_h$ is everywhere finite, and also provides a uniform local H\"older estimate for $\eps$-minimizers. For homogeneous potentials, the statement is due to Maderna \cite{Maderna2012}.

Here and below, we write $X \lesssim Y$ (resp. $X \gtrsim Y$) to indicate that there exists a constant $C > 0 $ such that $X \le CY$ (resp. $X \ge CY$); implied constants are allowed to depend on $V$ and on the masses $m_i$.

\begin{lemma}\label{lem:mane}
  For every $x , y \in (\RR^d)^N$,
  \begin{equation}\label{eq:ball}
    d_{h}(x,y)  \lesssim \|x - y\|_m^{\,1-a/2} \ +\ \|x - y\|_m.
  \end{equation}
  Moreover, for every $\eps \in [0,1]$ and every $\eps$-minimizer $\gamma:I \to (\RR^d)^N$,
  \begin{equation}\label{eq:holder}
    \left\|\gamma(s') - \gamma(s)\right\|_m  \lesssim \sqrt{s'-s}, \qquad \forall s,s'\in I\mbox{ with } 0<s'-s\le 2.
  \end{equation}
\end{lemma}

\begin{proof}

  The proof given in \cite{Maderna2012} can be adapted to our setting. However, we present an alternative argument which is inspired by the proof of Marchal's theorem (\Cref{lem:marchal} below) and uses a similar probabilistic argument. The idea is to join $x$ to $y$ by a broken line whose vertex is chosen at random, and to bound the expected $h$-length of the resulting random curve. As pointed out in \cite{Maderna2012}, the line segment from $x$ to $y$ itself may contain collisions or at least lie close to them; randomness is used in order to average out this scenario.

  \begin{figure}
    \centering
    \begin{tikzpicture}[x=1mm,y=1mm]
      \definecolor{diskgray}{HTML}{CCCCCC}
      \fill[diskgray] (0,0) circle[radius=22.489183];
      \draw[line width=0.246944mm,line cap=round,line join=round]
        (-22.45595,0.007) -- (5.32071,13.56358) -- (22.29573,0.12239);
      \fill (-22.45686,0.01857) circle[radius=0.375];
      \fill ( 22.45578,0.01894) circle[radius=0.375];
      \fill ( -0.01423,0.01894) circle[radius=0.375];
      \fill (  5.32071,13.56358) circle[radius=0.375];
      \node[anchor=base east] at (-22,-0.63661) {$x$};
      \node[anchor=base west] at ( 22,-0.63661) {$y$};
      \node[anchor=base]      at ( -0.99148,-2.6) {$z$};
      \node[anchor=base east] at ( -9.6,5.5)  {$\sigma^1$};
      \node[anchor=base west] at ( 15,5.7)      {$\sigma^2$};
      \node[anchor=south]     at ( 5.32071,13.5)  {$Z$};
    \end{tikzpicture}
    \caption{The random curve $\sigma$. Each component $Z_i$ of the vertex $Z$ is uniformly distributed on the ball of radius $\delta/2$ centered at $z_i$.}
    \label{fig:mane}
  \end{figure}

  We may assume that $x \ne y$. Write
  \[
    \delta : = \|x-y\|_{m} \qquad \text{and} \qquad z : = \frac{x+y}{2}.
  \]
  Let $Z = (Z_{1},\dots,Z_{N}) \in (\RR^{d})^{N}$ be a random vector whose coordinates $Z_{i} \in \RR^{d}$ are independent and uniformly distributed over the ball of radius $\delta/2$ centered at $z_i$. Let $\sigma$ denote the concatenation of the line segments $\sigma^1:[0,1]\to(\RR^{d})^{N}$ and $\sigma^2:[0,1]\to(\RR^{d})^{N}$, given by
  \[
    \sigma^1(t) : = (1-t)\,x + t\,Z \qquad \text{and} \qquad \sigma^2(t) : = (1-t)\,Z + t\,y,
  \]
  joining $x$ to $Z$ and $Z$ to $y$, respectively; see \Cref{fig:mane}.
  Then
  \begin{equation}
    \label{eq:mane-length}
    \begin{split}
      \Len_h(\sigma^1)  
      & =  \int_{0}^{1}\sqrt{2U\big(\sigma^1(t)\big)}\,\|\dot\sigma^1(t)\|_m\,dt \\
      & =  \int_{0}^{1}\sqrt{2U\big(\sigma^1(t)\big)}\,\|Z-x\|_{m}\,dt \\
      & \lesssim \delta\int_{0}^{1}\left(\sum_{i<j}f_{ij}\big(\|(1-t)(x_i - x_j) + t\,(Z_i - Z_j)\|_{\RR^d}\big)\right)^{1/2}dt\\
      & \lesssim \delta + \delta \int_0^1\ \sum_{i<j}\left\|(1-t)(x_{i}-x_{j}) + t\,(Z_i - Z_j)\right\|_{\RR^d}^{-a/2}\,dt,
    \end{split}
  \end{equation}
  where in the last inequality we used \eqref{eq:f-bounds}.
  We claim that, for every $i<j$ and $t\in(0,1]$,
  \begin{equation}\label{eq:maximized}
    v \mapsto \EE\left[\|v + t(Z_i - Z_j)\|_{\RR^d}^{-a/2}\right] \quad \text{attains its maximum at $v = -t(z_i - z_j)$. }
  \end{equation}
  Fix $i<j$ and $t \in (0,1]$, and set
  \[
    W : = t(Z_i - z_i) - t(Z_j - z_j).
  \]
  For each $w \in \RR^d$ and $r>0$, we denote by $B_r(w)$ the open ball of radius $r$ centered at $w$.
  The random vector $W$ has probability density
  \begin{align*}
    p(w)=\frac{\mathrm{vol}(B_{\delta/2}(0)\cap B_{\delta/2}(t^{-1}w))}{t^d\,\mathrm{vol}(B_{\delta/2}(0))^2}, \qquad w \in \RR^d,
  \end{align*}
  which is a nonincreasing function of $\|w\|_{\RR^d}$. Since 
  $\sup p|_{B_r(w) \setminus B_r(0)}
  \leq 
  \inf p|_{B_r(0)}$ and $\mathrm{vol}(B_r(w))=\mathrm{vol}(B_r(0))$,
  we have
  \begin{equation}
    \PP\big(W \in B_r(w)\big) \le \PP\big(W \in B_r(0)\big) .
  \end{equation}
  Applying this inequality with $w = -v - t(z_i - z_j)$ and $r = \lambda^{-2/a}$, together with the layer-cake formula, we get
  \begin{align*}
    \EE\big[\|v + t(Z_i - Z_j)\|_{\RR^d}^{-a/2}\big]
    &= 
    \EE\big[\|W + v + t(z_i - z_j)\|_{\RR^d}^{-a/2}\big]\\
    &=  \int_0^\infty \PP\big(\|W + v + t(z_i - z_j)\|_{\RR^d} < \lambda^{-2/a}\big)\,d\lambda\\
    &\le \int_0^\infty \PP\big(\|W\|_{\RR^d} < \lambda^{-2/a}\big)\,d\lambda
    = \EE\big[\|W\|_{\RR^d}^{-a/2}\big]\\
    &= \EE\big[\|-t(z_i - z_j) + t(Z_i - Z_j)\|_{\RR^d}^{-a/2}\big],
  \end{align*}
  which proves \eqref{eq:maximized}.
  By \eqref{eq:maximized}, for every $i<j$ and every $t \in (0,1]$,
  \begin{align*}
    \EE\left\|(1-t)(x_{i}-x_{j}) + t\,(Z_i - Z_j)\right\|_{\RR^d}^{-a/2}
    \le t^{-a/2}\,\EE\left\|(Z_i-Z_j) - (z_i - z_j)\right\|_{\RR^d}^{-a/2}  \lesssim (t\delta)^{-a/2}.
  \end{align*}
  Hence, taking expectations in \eqref{eq:mane-length} and using Fubini's theorem, we obtain
  \[
    \EE\big[\Len_h(\sigma^1)\big]  \lesssim \delta + \delta^{1-a/2}\int_{0}^{1}t^{-a/2}\,dt  \lesssim \delta + \delta^{\,1-a/2}
  \]
  since $a < 2$.
  The same argument shows that $\EE\left[\Len_h(\sigma^2)\right]\lesssim \delta + \delta^{\,1-a/2}$, and therefore
  \[
    d_h(x,y) \le \EE[\Len_h(\sigma)]\lesssim \delta + \delta^{\,1-a/2},
  \]
  which is \eqref{eq:ball}.

  Let $\eps \in [0,1]$, let $\gamma:I \to (\RR^d)^N$ an $\eps$-minimizer, and fix $s,s' \in I$ with $0 < s'-s \le 2$. By \eqref{eq:maupertuis}, positivity of $U_\eps$, and Cauchy--Schwarz,
  \begin{equation}\label{eq:holder-CS}
    d_{h_\eps}(\gamma(s),\gamma(s')) = \cA_\eps\!\big(\gamma\vert_{[s,s']}\big) \ge \frac12\int_{s}^{s'}\left\|\dot\gamma(t)\right\|_m^2dt \ge \frac{\left\|\gamma(s')-\gamma(s)\right\|_m^2}{2(s'-s)},
  \end{equation}
  so since $s'-s \le 2$,
  \begin{equation}\label{dhepsgamma1}
    d_{h_\eps}(\gamma(s),\gamma(s')) \gtrsim \|\gamma(s') - \gamma(s)\|_m^2.
  \end{equation}
  On the other hand, $d_{h_\eps} \le d_h$  by \eqref{eq:reg}, so by \eqref{eq:ball} applied to $x = \gamma(s)$ and $y = \gamma(s')$,
  \begin{equation}\label{dhepsgamma2}
    d_{h_\eps}(\gamma(s),\gamma(s'))
    \le 
    d_{h}(\gamma(s),\gamma(s'))
    \lesssim 
    \left\|\gamma(s')-\gamma(s)\right\|_m^{\,1-a/2} + \left\|\gamma(s')-\gamma(s)\right\|_m.
  \end{equation}
  Comparing \eqref{dhepsgamma1} and \eqref{dhepsgamma2} gives $\left\|\gamma(s')-\gamma(s)\right\|_m \lesssim 1$, hence $d_{h_\eps}(\gamma(s),\gamma(s')) \lesssim 1$ by \eqref{dhepsgamma2}, and then $\left\|\gamma(s')-\gamma(s)\right\|_m^2 \lesssim s'-s$ by \eqref{eq:holder-CS}, which is \eqref{eq:holder}.
\end{proof}

By \Cref{lem:mane}, there exists a constant $C_1 > 0$ such that for every $\eps \in [0,1]$ and every $\eps$-minimizer $\gamma:I \to (\RR^d)^N$ defined on an interval $I \ni 0$,
\begin{equation}\label{eq:holder-pairs}
  \left\|(\gamma_i - \gamma_j)(s) - (\gamma_i - \gamma_j)(0)\right\|_{\RR^d} \le C_1\,\sqrt{|s|}, \qquad \forall 1 \le i<j \le N,\quad \forall s \in I\cap [-2,2].
\end{equation}

Recall that our aim is to show that an entire minimizer is also an entire $\eps$-minimizer if $\eps>0$ is sufficiently small.
This could fail if we had an entire minimizer going through a collision, since the potentials $V_\eps$ and $V$ would not coincide in a neighborhood of it. We shall show that, indeed, entire minimizers always avoid collisions.
For the Newtonian potential, the fact that free-time minimizers avoid collisions in the interior of their domains is known as \emph{Marchal's theorem} \cite{Marchal2002,Chenciner2002,FerrarioTerracini2004,Montgomery2018}.
We shall establish a variant of Marchal's theorem which gives a lower bound on the distance of $\eps$-minimizers from collisions, uniformly over $\eps$.

\begin{lemma}[Uniform Marchal's theorem]\label{lem:marchal}
  For every $T \in (0,1]$ there exists $R \in (0,1)$ such that for every $\eps \in [0,R]$ and every $\eps$-minimizer $\gamma : [-T,T] \to (\RR^d)^N$,
  \[
    \left\|\gamma_i(0) - \gamma_j(0)\right\|_{\RR^d}  > R, \qquad \forall 1 \le i<j \le N.
  \]
\end{lemma}

\begin{proof}
  Fix $T \in (0,1]$ and $1 \le i<j \le N$. By \eqref{eq:f-singular}, the difference $f_{ij}(r/4) - f_{ij}(r/2)$ is unbounded as $r \to 0$. We may therefore fix $R_{ij} \in (0,1)$ such that
  \begin{equation}\label{eq:R-choice}
    R_{ij}^2 \le \frac{T}{C_2}
    \qquad\text{and}\qquad
    f_{ij}\!\left(R_{ij}/4\right) - f_{ij}\!\left(R_{ij}/2\right) \ >\ \frac{C_2}{T} ,
  \end{equation}
  where $C_2 > 0$ is a constant whose value will be determined below. 
  We will show that if
  \begin{equation}\label{eq:epschoice}
    \eps \le \frac{R_{ij}}{4}
  \end{equation}
  then every $\eps$-minimizer $\gamma : [-T,T] \to (\RR^d)^N$ satisfies $\|\gamma_i(0) - \gamma_j(0)\|_{\RR^d} > R_{ij}/8$; the lemma thus holds with $R : = \min_{i<j}R_{ij}/8$.

  Suppose to the contrary that there exists $\eps \le R_{ij}/4$ and an $\eps$-minimizer $\gamma:[-T,T]\to(\RR^{d})^N$ such that
  \begin{equation}\label{eq:gammaassump}
    \|\gamma_i(0) - \gamma_j(0)\|_{\RR^d} \le R_{ij}/8.
  \end{equation}
  Following Marchal \cite{Marchal2002,Chenciner2002}, we will construct a random curve $\eta$ with the same endpoints as $\gamma$ such that
  \begin{equation}\label{eq:EEAgammav}
    \EE\!\left[\cA_\eps(\eta)\right] < \cA_\eps(\gamma);
  \end{equation}
  in particular, $\cA_\eps(\eta) < \cA_\eps(\gamma)$ with positive probability, contradicting the minimality of $\gamma$.

  Let $\zeta : \RR \to [0,1]$ be a smooth function supported in $(-T,T)$ such that $\zeta\vert_{[-T/2,T/2]} \equiv 1$ and $|\dot\zeta| \le 4/T$, so that 
  \begin{equation}\label{eq:zeta-energy}
    \int_{-T}^{T}\dot\zeta(t)^2\,dt  \le  \frac{16}{T} .
  \end{equation}
  Let $\theta$ be a random vector distributed uniformly over the unit sphere $S^{d-1}\subseteq\RR^d$, and define  $\eta:[-T,T] \to (\RR^d)^N$ by
  \[
    \eta_{k}(t)
    : = 
    \begin{cases}
      \gamma_{k}(t), & \mbox{if }k \neq j \\
      \gamma_{j}(t) + \zeta(t)\cdot R_{ij}\theta, & \mbox{if }k = j .
    \end{cases}
  \]
  Since $\EE\,\theta = 0$ and $\|\theta\|_{\RR^d} = 1$, 
  \[
    \EE\,\|\dot\gamma_j + \dot\zeta\cdot R_{ij}\theta\|_{\RR^d}^2 = \|\dot\gamma_j\|_{\RR^d}^2 + R_{ij}^2\dot\zeta^2.
  \]
  Since moreover $\cA_\eps(\gamma) < \infty$ by \eqref{eq:ball}, the definition of $\cA_\eps$ and Fubini's theorem give (setting $f^{\eps}_{lk} : = f^{\eps}_{kl}$ for $k < l$):
  \begin{equation}\label{eq:marchal-average}
    \begin{split}
      & \EE\left[\cA_\eps(\eta)\right] - \cA_\eps(\gamma)\\
      & \quad  =  \EE\int_{-T}^{T}\left(\tfrac12\left\|\dot\eta\right\|_m^2 - V_\eps(\eta)\right)dt \ - \ \int_{-T}^{T}\left(\tfrac12\left\|\dot\gamma\right\|_m^2 - V_\eps(\gamma)\right)dt \\
      & \quad  =  \EE\int_{-T}^{T}\bigg[\frac{m_j}{2}\Big(\|\dot\gamma_j + \dot\zeta\cdot R_{ij}\theta\|_{\RR^d}^2 - \|\dot\gamma_j\|_{\RR^d}^2\Big) \\
      & \qquad  + \sum_{k \neq j}\Big(f^{\eps}_{kj}\!\left(\left\|\gamma_{k} - \gamma_{j} - \zeta\cdot R_{ij}\theta\right\|_{\RR^d}\right) - f^{\eps}_{kj}\!\left(\left\|\gamma_{k} - \gamma_{j}\right\|_{\RR^d}\right)\Big)\bigg]dt \\
      & \quad =  \frac{m_j}{2}\,R_{ij}^2\int_{-T}^{T}\dot\zeta(t)^2\,dt \\
      & \qquad  + \sum_{k \neq j}\int_{-T}^{T}\bigg(\EE f^{\eps}_{kj}\!\left(\left\|\gamma_{k}(t) - \gamma_{j}(t) - \zeta(t)\cdot R_{ij}\theta\right\|_{\RR^d}\right) - f^{\eps}_{kj}\!\left(\left\|\gamma_{k}(t) - \gamma_{j}(t)\right\|_{\RR^d}\right)\bigg)dt .
    \end{split}
  \end{equation}
  The function $f^{\eps}_{kj}(\|\cdot\|_{\RR^d})$ is superharmonic on $\RR^d$, so its average over a sphere is at most its value at the center of the sphere. Thus each of the integrands in the last line of \eqref{eq:marchal-average} is nonpositive, and we may discard all of them except for $k = i$ and $|t| \le t_0$, where
  \[
    t_{0} : = \frac{R_{ij}^{2}}{64\,C_1^2} \overset{\eqref{eq:R-choice}}{\le} \frac{T}{64\,C_1^2C_2} \le \frac{T}{2},
  \]
  where the last inequality holds true provided that $C_2 \ge 1/(32\,C_1^2)$. Hence
  \[
    \zeta\vert_{[-t_0,t_0]} \equiv 1.
  \]
  Moreover, by \eqref{eq:holder-pairs} and \eqref{eq:gammaassump}, the function $w : = \gamma_i - \gamma_j$ satisfies
  \[
    \|w(t)\|_{\RR^d} \le \|w(0)\|_{\RR^d} + C_1\sqrt{|t|} \le R_{ij}/4,
  \]
  and hence
  \[
    \|w(t) - R_{ij}\theta\|_{\RR^d}  \ge R_{ij}/2, \qquad  \forall |t| \le t_0.
  \]
  Since $f^{\eps}_{ij}$ is nonincreasing and coincides with $f_{ij}$ on $[R_{ij}/4,\infty)$ by \eqref{eq:reg} and \eqref{eq:epschoice}, we obtain
  \begin{align*}
    \EE\left[\cA_\eps(\eta)\right] - \cA_\eps(\gamma)
    \ & \le \ \frac{m_j}{2}\,R_{ij}^2\int_{-T}^{T}\dot\zeta(t)^2\,dt   \\ 
    & \quad +\int_{-t_0}^{t_{0}}\bigg(\EE f^{\eps}_{ij}\!\left(\left\|w(t) - R_{ij}\theta\right\|_{\RR^d}\right) - f^{\eps}_{ij}\!\left(\left\|w(t)\right\|_{\RR^d}\right)\bigg)dt \\
    \ & \le \ \frac{m_j}{2}\,R_{ij}^2\int_{-T}^{T}\dot\zeta(t)^2\,dt \ - \ 2t_0\bigg(f_{ij}\!\left(R_{ij}/4\right) - f_{ij}\!\left(R_{ij}/2\right)\bigg) \\
    \ & \le \ R_{ij}^2\left(\frac{8m_j}{T} \ - \ \frac{f_{ij}\!\left(R_{ij}/4\right) - f_{ij}\!\left(R_{ij}/2\right)}{32\,C_1^2}\right) ,
  \end{align*}
  where the last inequality uses \eqref{eq:zeta-energy} and the definition of $t_0$. 
  Hence, by \eqref{eq:R-choice},
  \[
    \EE\left[\cA_\eps(\eta)\right] - \cA_\eps(\gamma) < \frac{R_{ij}^2}{T} \cdot \left(8m_j - \frac{C_2}{32\,C_1^2}\right) \le 0
  \]
  provided that $C_2 \ge 256\,m_jC_1^2$. Thus with the choice
  \[
    C_2 : = \max\left\{\frac{1}{32\,C_1^2},\ 256\,m_jC_1^2\right\}
  \]
  inequality \eqref{eq:EEAgammav} holds, and the lemma is proved.
\end{proof}

\begin{lemma}\label{clear:main}
  Suppose that there exists an entire minimizer $\gamma$. Then, for all $\eps$ sufficiently small, $\gamma$ is an entire $\eps$-minimizer as well.
\end{lemma}

\begin{proof}
  For every $s \in \RR$ the curve $t \mapsto \gamma(s+t)$, $t \in [-1,1]$, is a minimizer. Hence \Cref{lem:marchal}, applied with $T = 1$ and $\eps = 0$, shows that
  \begin{equation}\label{eq:D-def}
    D : = \inf_{s\in \RR,\, i<j}\left\|\gamma_i(s)-\gamma_j(s)\right\|_{\RR^d} \ > \ 0 .
  \end{equation}
  Fix $\tau \in (0,1]$ such that $C_1\sqrt{\tau} \le D/2$, where $C_1$ is the constant from \eqref{eq:holder-pairs}.
  Let $R \in (0,1)$ be the constant given by \Cref{lem:marchal} with $T = \tau$, and set $\eps_0 : = \min\{R, D/2\}$.

  Fix $\eps \in (0,\eps_0]$ and $t < t'$. 
  The smooth metric $h_\eps$ is complete by \Cref{reg:hypotheses}, so there is a minimizing $h_\eps$-geodesic joining $\gamma(t)$ to $\gamma(t')$; let $\eta : [0,\ell] \to (\RR^d)^N$ be its parametrization by unit $h_\eps^*$-speed, which is an $\eps$-minimizer by \eqref{eq:length-vs-action} and \eqref{eq:maupertuis}. 
  We claim that
  \begin{equation}\label{eq:eta-away}
    \left\|\eta_i(s)-\eta_j(s)\right\|_{\RR^d} \ge \eps_0, \qquad \forall s \in [0,\ell],\ i<j.
  \end{equation}
  For $s \in [\tau,\ell-\tau]$ the curve $t \mapsto \eta(s+t)$, $t \in [-\tau,\tau]$, is an $\eps$-minimizer and $\eps \le \eps_0 \le R$, so \Cref{lem:marchal} gives $\left\|\eta_i(s)-\eta_j(s)\right\|_{\RR^d} > R \ge \eps_0$, which is \eqref{eq:eta-away} for such $s$. Suppose therefore that $s \in [0,\tau)$ (the case $s \in (\ell-\tau,\ell]$ is analogous). By \eqref{eq:holder-pairs}, since $\eta$ is an $\eps$-minimizer and $0 \le s < \tau \le 1$,
  \[
    \|(\eta_i - \eta_j)(s) - (\eta_i - \eta_j)(0)\|_{\RR^d} \le C_1\sqrt{\tau}.
  \]
  Since $\eta(0) = \gamma(t)$, it follows that
  \[
    \left\|\eta_i(s)-\eta_j(s)\right\|_{\RR^d} \ge \left\|\gamma_i(t)-\gamma_j(t)\right\|_{\RR^d} - C_1\sqrt{\tau} \ge D - \frac{D}{2} \ge \eps_0,
  \]
  and \eqref{eq:eta-away} is proved.
  By \eqref{eq:reg}, \eqref{eq:D-def}, \eqref{eq:eta-away} and the fact that $\eps_0 \le D$, the potentials $V_\eps$ and $V$ coincide both on the curve $\gamma$ and on the curve $\eta$ for all $\eps \le \eps_0$.
  Since $\gamma$ is a minimizer, it follows that $\cA_\eps(\gamma\vert_{[t,t']}) = \cA_0(\gamma\vert_{[t,t']}) \le \cA_0(\eta) = \cA_\eps(\eta)$. But $\eta$ is an $\eps$-minimizer, hence so is $\gamma\vert_{[t,t']}$. Since $t$ and $t'$ are arbitrary, $\gamma$ is an $\eps$-minimizer.
\end{proof}

\begin{proof}[Proof of \Cref{thm:mb}]
  Assume by contradiction that there exists an entire minimizer $\gamma : \RR \to (\RR^d)^N$.
  By \Cref{clear:main}, we may fix $\eps \in (0,1]$ such that $\gamma$ is an entire $\eps$-minimizer. By \eqref{eq:calibrated} and \eqref{eq:unit-speed},
  \begin{equation}\label{eq:line-length}
    d_{h_\eps}(\gamma(t),\gamma(t'))  =  \cA_\eps\big(\gamma\vert_{[t,t']}\big)  =  \Len_{h_\eps}\!\big(\gamma\vert_{[t,t']}\big)  =  \int_t^{t'}2U_\eps(\gamma(s))\,ds, \qquad \forall t < t'.
  \end{equation}
  Since $h_\eps$ is complete by \Cref{reg:hypotheses} and $U_\eps$ is positive, $d_{h_\eps}(\gamma(0),\gamma(t')) \to \infty$ as $t' \to \infty$. Indeed, if $d_{h_\eps}(\gamma(0),\gamma(t'))$ were bounded then $\gamma([0,\infty))$ would lie in a compact set, on which $U_\eps$ is bounded below by a positive constant, and the integral in \eqref{eq:line-length} would diverge. Similarly, $d_{h_\eps}(\gamma(t),\gamma(0)) \to \infty$ as $t \to -\infty$.  Hence the reparametrization of $\gamma$ by $h_\eps$-arclength is a line for the metric $h_\eps$, so \Cref{reg:hypotheses} and \Cref{thm:main} imply that $U_\eps$ is constant, which is a contradiction since the functions $f_{ij}$ are strictly decreasing and $f_{ij}^\eps(r) = f_{ij}(r)$ for large $r$.
\end{proof}

\section[Proof of Theorem 2]{Proof of \Cref{thm:main}}\label{sec:proof}

This section follows closely the proof of the Cheeger--Gromoll splitting theorem \cite{CheegerGromoll1971,EschenburgHeintze1984,Petersen2016}. For the reader's convenience, we will provide a detailed proof.

Let $(M,g)$ be a connected Riemannian manifold of dimension $n \ge 2$ and let $U : M \to (0,\infty)$ be a smooth function. Throughout this section we assume that \eqref{eq:DeltaVcondition} holds and that the conformal metrics $h = Ug$ and $h^* = U^{-1}g$ are both complete.
The operators $\Delta$, $\nabla$, $\Hess$ and $\Ric$ are taken with respect to the metric $g$.

\begin{lemma}[{cf.\ \cite[Lemma 1]{EschenburgHeintze1984}}]
  Let $W \subseteq M$ be an open set and let $u : W \to \RR$ be a smooth solution to 
  \begin{equation}\label{eq:HJ}
    \|du\|_g^2 = U .
  \end{equation}
  Then
  \begin{equation}\label{eq:riccati}
    (d\Delta u)(\nabla u) \le -\frac{(\Delta u)^2}{n} ,
  \end{equation}
  with equality at $q \in W$ if and only if
  \begin{equation}\label{eq:riccati-equality}
    \begin{cases}
      \tfrac12\,\Delta U(q) = \Ric(\nabla u\vert_q,\nabla u\vert_q),\vspace{5pt}\\
      \Hess u\vert_q = \displaystyle\frac{\Delta u(q)}{n}\,g\vert_q .
    \end{cases}
  \end{equation}
\end{lemma}

\begin{proof}
  By the Cauchy--Schwarz inequality,
  \begin{equation}\label{eq:Hessinequ}
    \|\!\Hess u\|_g^2 \ge \frac{(\Delta u)^2}{n}, 
  \end{equation}
  with equality at $q \in W$ if and only if $\Hess u\vert_q = \frac{\Delta u(q)}{n}\,g\vert_q$.
  The Bochner--Weitzenb\"ock formula
  \[
    \tfrac12\Delta\|d u\|_g^2 = (d\Delta u)(\nabla u) + \|\!\Hess u\|_g^2 + \Ric(\nabla u,\nabla u),
  \]
  together with \eqref{eq:HJ}, gives
  \begin{align*}
    (d\Delta u)(\nabla u) 
    & = \tfrac12\,\Delta U - \Ric(\nabla u,\nabla u) - \|\!\Hess u\|_g^2 \\
    & = \tfrac1{U}\left(\tfrac12\,(\Delta U)\,g - U\Ric\right)(\nabla u,\nabla u) - \|\!\Hess u\|_g^2,
  \end{align*}
  whence by \eqref{eq:DeltaVcondition}
  \begin{equation}\label{eq:second-order}
    (d\Delta u)(\nabla u) \le - \|\!\Hess u\|_g^2.
  \end{equation}
  Inequalities \eqref{eq:Hessinequ} and \eqref{eq:second-order} combine to give \eqref{eq:riccati}.
  Equality holds in \eqref{eq:riccati} at $q \in W$ if and only if equality holds at $q$ in \eqref{eq:Hessinequ} and $\tfrac12\,(\Delta U)\,\|\nabla u\|_g^2 - U \cdot \Ric(\nabla u,\nabla u) = 0$ at $q$, which is \eqref{eq:riccati-equality}.
\end{proof}

\begin{lemma}\label{lem:harmonic}
  If there exists a harmonic function $u:M \to \RR$ satisfying \eqref{eq:HJ}, then $U$ is constant.
\end{lemma}

\begin{proof}
  Suppose that $u : M \to \RR$ is harmonic and solves \eqref{eq:HJ}. Then both sides of \eqref{eq:riccati} vanish, so \eqref{eq:riccati-equality} holds everywhere and $\Hess u \equiv 0$. Hence, by \eqref{eq:HJ},
  \[
    dU = d\big(\|du\|_g^2\big) = 2\Hess u\,(\nabla u,\cdot\,) \equiv 0 .
  \]
  Since $M$ is connected, it follows that $U$ is constant.
\end{proof}

For $q \in M$, let $\Cut_h(q) \subseteq M$ denote the cut locus of $q$ with respect to the metric $h$, and write
\[
  u_q := d_h(q,\cdot\,).
\]

\begin{lemma}[Laplacian comparison, cf.\ {\cite[Lemma 2]{CheegerGromoll1971}}]\label{lem:comparison}
  Let $\gamma : [0,T] \to M$ be a minimizing $h$-geodesic such that $\gamma(t) \notin \Cut_h(\gamma(0))$ for all $t\in(0,T]$. Then
  \[
    \Delta u_{\gamma(0)}(\gamma(t)) \le \frac{n}{\Len_{h^*}\!\big(\gamma\vert_{[0,t]}\big)} , \qquad \forall t\in(0,T].
  \]
\end{lemma}

\begin{proof}
  We set $q:=\gamma(0)$.
  Away from $\{q\}\cup\Cut_h(q)$, the function $u_q = d_h(q,\cdot)$ is smooth with 
  \begin{align}\label{eq:h_norm_duq}
    \|du_q\|_h \equiv 1,
  \end{align}
  and therefore solves \eqref{eq:HJ}.
  Since $\gamma$ is a minimizing $h$-geodesic issuing from $q$, we have $\gamma(t) \ne q$ for all $t\in(0,T]$.
  This, together with the assumption that $\gamma(t)\not\in\Cut_h(q)$ for all $t\in(0,T]$, implies that $u_q$ is smooth near $\gamma(t)$ for all $t\in(0,T]$, and $\dot\gamma$ is positively proportional to the $h$-gradient of $u_q$. Since $h$ is conformal to $g$, it follows that $\dot\gamma$ is a positive multiple of $\nabla u_q\circ\gamma$. 
  Let $\sigma:[0,T']\to M$ be the reparametrization of $\gamma$ with unit $h^*$-speed, i.e.\ $\|\dot\sigma\|_{g}^2=U(\sigma)$. Now, $\dot\sigma$ is a positive multiple of $\nabla u_q\circ\sigma$ and, by \eqref{eq:h_norm_duq}, both have $g$-norm 
  \[\|\dot\sigma\|_g=\|\nabla u_q(\sigma)\|_g=\sqrt{U(\sigma)}.\]
  Thus
  \[
    \dot\sigma = \nabla u_q\circ\sigma,
  \]
  whence by \eqref{eq:riccati} the function $\psi : = \Delta u_q\circ\sigma$ satisfies
  \[
    \dot\psi \le -\psi^2/n \le 0,
  \]
  and we need to prove that $\psi(t)\le n/t$ for all $t \in (0,T']$. If $\psi(t) \le 0$ there is nothing to prove; otherwise, since $\psi$ is nonincreasing, necessarily $\psi > 0$ on $(0,t]$, so $\tfrac{d}{ds}\psi(s)^{-1} \ge 1/n$ there. Integrating from $s$ to $t$ and letting $s \to 0^+$ gives $\psi(t)^{-1} \ge t/n$.
\end{proof}

Assume from now on that there exists a line $\gamma : \RR \to M$ for the metric $h$. For $t \in \RR$, set
\[
  b_t : = d_h(\gamma(t),\cdot\,) - d_h(\gamma(t),\gamma(0)) .
\]
By the triangle inequality, each $b_t$ is $1$-Lipschitz with respect to the metric $h$, and $t \mapsto b_t(q)$ is nonincreasing on $[0,\infty)$, nondecreasing on $(-\infty,0]$, and bounded below by $-d_h(q,\gamma(0))$ for all $q \in M$. Hence the \emph{Busemann functions}
\[
  b_\pm : M \to \RR, \qquad b_\pm : = \lim_{t\to\pm\infty}b_t
\]
are well-defined and $1$-Lipschitz with respect to the metric $h$.
For any $s \in \RR$ we have $b_t(\gamma(s)) = -s$ if $t > \max\{s,0\}$, while $b_t(\gamma(s)) = s$ if $t < \min\{s,0\}$; letting $t \to \pm\infty$,
\begin{equation}\label{eq:horo-on-gamma}
  b_+(\gamma(s)) = -s = -\,b_-(\gamma(s)) , \qquad \forall s \in \RR .
\end{equation}
Moreover, since $\gamma$ is a line, the triangle inequality gives $b_t + b_{t'} \ge 0$ for all $t < 0 < t'$, whence
\begin{equation}\label{eq:horo-sum}
  b_+(q) + b_-(q) \ge 0, \qquad\forall q\in M .
\end{equation}

\begin{lemma}[{cf.\ \cite[Prop.~7.3.8]{Petersen2016}, \cite[\S2]{EschenburgHeintze1984}, \cite[Thm.~2]{CheegerGromoll1971}}]\label{lem:asymptote}
  For every $q \in M$ there is a unit-speed minimizing $h$-geodesic $\sigma : [0,\infty) \to M$ with initial point $\sigma(0) = q$ and such that
  \begin{equation}\label{eq:unit-rate}
    b_+(\sigma(t)) = b_+(q) - d_h(q,\sigma(t)) , \qquad \forall t \ge 0 .
  \end{equation}
  The same statement holds for $b_-$.
\end{lemma}

\begin{proof}
  Since $h$ is complete, for every $t > 0$ such that $\gamma(t)\neq q$ there exists a unit-speed minimizing $h$-geodesic $\sigma_t$ such that $\sigma_t(0) = q$ and $\sigma_t(\ell_t) = \gamma(t)$. Notice that 
  \[
    \ell_t = d_h(q,\gamma(t)) \ge\ t - d_h(q,\gamma(0)).
  \]
  Since each $\dot\sigma_t(0)$ is a vector in $T_qM$ of norm $\|\dot\sigma_t(0)\|_h=1$, there exists a sequence $t_k \to \infty$ such that $\dot\sigma_{t_k}(0)\to v\in T_qM$ with $\|v\|_h=1$. Since $h$ is complete, the unit-speed $h$-geodesic $\sigma : [0,\infty) \to M$ with $\sigma(0) = q$ and $\dot\sigma(0) = v$ is defined on the whole half-line, and $\sigma_{t_k} \to \sigma$ locally uniformly on $[0,\infty)$ by continuous dependence of geodesics on their initial velocity.
  For every $k \ge 1$ and every $r \in [0,d_h(q,\gamma(t_k))]$,
  \begin{align*}
    b_{t_k}(\sigma_{t_k}(r)) &= d_h(\gamma(t_k),\sigma_{t_k}(r)) - d_h(\gamma(t_k),\gamma(0))\\
    & = d_h(\gamma(t_k),q) - r - d_h(\gamma(t_k),\gamma(0)) = b_{t_k}(q) - r. 
  \end{align*}
  Letting $k \to \infty$, and using that $b_{t_k}$ is $1$-Lipschitz with respect to the metric $h$ and $\sigma_{t_k}(r) \to \sigma(r)$, we infer
  \begin{align*}
    b_+(\sigma(r)) &= \lim_{k\to\infty}b_{t_k}(\sigma(r)) = \lim_{k\to\infty}b_{t_k}(q) - r = b_+(q) - r,\\
    \quad d_h(q,\sigma(r)) &= \lim_{k\to\infty}d_h(q,\sigma_{t_k}(r)) = r, 
  \end{align*}
  for all $r \ge 0$, and \eqref{eq:unit-rate} follows. Finally, since $d_h(q,\sigma(r)) = r = \Len_{h}\!\big(\sigma\vert_{[0,r]}\big) $ for all $r \ge 0$, the $h$-geodesic $\sigma$ is minimizing.
  The proof for $b_-$ is similar.
\end{proof}

\begin{lemma}\label{lem:horo-HJ}
  At every point $q\in M$ of differentiability of $b_\pm$, we have
  \[
    \|db_\pm(q)\|_g^2 = U(q).
  \]
\end{lemma}

\begin{proof}
  Let $q \in M$ be a point of differentiability of $b_+$ and let $\sigma$ be as in \Cref{lem:asymptote}. Since $b_+$ is differentiable at $q$, the right derivative of $t \mapsto b_+(\sigma(t))$ at $t = 0$ is $db_+(q)\dot\sigma(0)$, and by \eqref{eq:unit-rate} it equals $-1 = -\|\dot\sigma(0)\|_h$. Hence $\|db_+(q)\|_h \ge 1$. The reverse inequality holds since $b_\pm$ are $1$-Lipschitz with respect to $h$.
\end{proof}

Following Calabi \cite[Def.~1]{Calabi1958}, a continuous function $u$ on an open set $W\subseteq M$ is said to satisfy $\Delta u \ge 0$ \emph{in the support sense} at a point $q \in W$ if for every $\eps > 0$ there exists a $C^2$ function $v_\eps$ defined on a neighborhood of $q$ such that
$v_\eps \le u$, $v_\eps(q) = u(q)$, and $\Delta v_\eps(q) \ge -\eps$.
We say that $\Delta u \ge 0$ in the support sense on $W$ if this holds at every point of $W$. We say that $\Delta u \le 0$ in the support sense if $\Delta(-u) \ge 0$ in the support sense.

\begin{proposition}[{cf.\ \cite[Thm.~1]{CheegerGromoll1971}, \cite[\S5]{EschenburgHeintze1984}}]\label{prop:superharmonic}
  $\Delta b_\pm \le 0$ in the support sense on $M$.
\end{proposition}

\begin{proof}
  Fix $q \in M$ and let $\sigma$ be as in \Cref{lem:asymptote}. Since $\sigma$ is a minimizing $h$-geodesic, $\sigma(s) \notin \Cut_h(\sigma(t))$ for all $0 \le s < t$. For $t > 0$ set
  \[
    v_t : = b_+(\sigma(t)) + d_h(\sigma(t),\cdot\,) = b_+(\sigma(t)) + u_{\sigma(t)}.
  \]
  Then $v_t$ is smooth near $q$ since $q \notin \Cut_h(\sigma(t))$, and the triangle inequality implies that $b_+ \le v_t$, with equality at $q$ by \eqref{eq:unit-rate}. Applying \Cref{lem:comparison} to the curve $\sigma\vert_{[0,t]}$ traversed backwards, we see that
  \[
    \Delta v_t(q) = \Delta u_{\sigma(t)}(q) \le \frac{n}{\Len_{h^*}\!\big(\sigma\vert_{[0,t]}\big)} .
  \]
  Since $\sigma:[0,\infty)\to M$ is a unit-speed minimizing $h$-geodesic, it must leave every compact subset of $M$. This, together with the completeness of the metric $h^*$, implies that $\Len_{h^*}\!\big(\sigma\vert_{[0,t]}\big) \to \infty$ as $t \to \infty$. We conclude that $\Delta b_+ \le 0$ in the support sense. The argument for $b_-$ is analogous.
\end{proof}

\begin{proof}[Proof of \Cref{thm:main}]
  Define $u:M \to \RR$ by $u : = b_+ + b_-$.
  Then $u \ge 0$ by \eqref{eq:horo-sum}, and $u(\gamma(0))=0$ by \eqref{eq:horo-on-gamma}. Hence $u$ attains a minimum at $\gamma(0)$. By \Cref{prop:superharmonic}, $\Delta b_\pm \le 0$ in the support sense, and therefore $\Delta u \le 0$ in the support sense as well. By Calabi's extension of E.~Hopf's maximum principle {\cite[Thm.~1]{Calabi1958}}, the function $u$ must vanish identically, i.e.
  \[
    b_+(q) =  -b_-(q),\qquad \forall q\in M .
  \]
  Therefore, we have both $\Delta b_+\leq0$ and $\Delta b_+\geq0$  in the support sense. 
  It follows that $b_+$ is harmonic (cf.\ \cite[p.~120]{CheegerGromoll1971}). Indeed, fix $x \in M$ and a small enough closed $g$-geodesic ball $B \subset M$ centered at $x$, and let $f:B\to\RR$ be the harmonic function such that $f|_{\partial B}=b_+|_{\partial B}$. As both $\Delta(f - b_+) \le 0$ and $\Delta(f - b_+) \ge 0$ in the support sense and $f|_{\partial B}-b_+|_{\partial B}\equiv0$, the above maximum principle gives $(f-b_+)|_{B}\equiv0$. Therefore $b_+|_B = f$ is harmonic, and in particular smooth. Since $x$ is arbitrary, $b_+$ is harmonic on all of $M$. By \Cref{lem:horo-HJ} the function $b_+$ is a solution to \eqref{eq:HJ}, and so by \Cref{lem:harmonic} the function $U$ is identically equal to some constant $c > 0$.

  Thus $h = c\,g$, so that $(M,g)$ itself is complete, and assumption \eqref{eq:DeltaVcondition} reads $\Ric_g \ge 0$. The rest of the argument is the classical one, cf.\ \cite[Remark 2]{EschenburgHeintze1984}, \cite[p.~121]{CheegerGromoll1971} or \cite[Thm.~7.3.5]{Petersen2016}.
  The function $\tau : = c^{-1/2}\,b_+$ is harmonic, hence smooth, and has unit gradient by \Cref{lem:horo-HJ}. As in the proof of \Cref{lem:harmonic}, both sides of \eqref{eq:riccati} vanish, so that \eqref{eq:riccati-equality} gives $\Hess\tau \equiv 0$, i.e., the gradient $\nabla\tau$ is a parallel unit vector field, and in particular a Killing vector field. Since $(M,g)$ is complete and $\|\nabla\tau\|_g \equiv 1$, the function $\tau : M \to \RR$ is surjective and the flow $\Phi_t$ of $\nabla\tau$ is complete. Let $M' : = \tau^{-1}(0)$, endowed with the metric $g'$ induced by $g$. Since $\tau(\Phi_t(x)) = \tau(x) + t$, the map $\iota(x,t):=\Phi_t(x)$ is a diffeomorphism from $M'\times \RR$ to $M$. As $M'$ is orthogonal to $\nabla\tau$ and each $\Phi_t$ is an isometry, $\iota^*g = g' + dt^2$. Finally $(M',g')$ is complete because $(M,g)$ is, and $\Ric_{g'} \ge 0$ because $\Ric \ge 0$.
\end{proof}

\bibliography{biblio}

\vspace{10pt}

\end{document}